\magnification=\magstep1
\vbadness=10000
\hbadness=10000
\tolerance=10000

\def\R{{\bf R}}   

\nopagenumbers

\proclaim What is a vertex algebra?.

Arbeitstagung talk 1997

Richard E. Borcherds, 
\footnote{$^*$}{ Supported by a Royal Society
professorship.}

D.P.M.M.S.,
16 Mill Lane, 
Cambridge, 
CB2 1SB,
England.

reb@dpmms.cam.ac.uk, http://www.dpmms.cam.ac.uk/$\tilde{~}$reb

\bigskip

The answer to the question in the title is that a vertex algebra is
really a sort of commutative ring. I will try to explain this in the
rest of the talk, and show how to use this to generalize the idea of a
vertex algebra to higher dimensions. The picture to keep in mind is
that a commutative ring should be thought of as somehow related to
quantum field theories in 0 dimensions, and vertex algebras are
related in the same way to 1 dimensional quantum field theories, and
we want to find out what corresponds to higher dimensional field
theories.  This talk is an exposition of the paper q-alg/9706008,
which contains (some of) the missing details. There is also
some overlap with unpublished notes of Soibelman, which have just
appeared on the q-alg preprint server as q-alg/9709030.

The relation of vertex algebras to commutative rings is obscured by 
the rather bad notation generally used for vertex algebras. Recall that 
for any element $v$ of a vertex algebra $V$ we have a vertex operator
denoted by $V(v,z)$ taking $V$ to the Laurent power series in $V$. 
I am going to change notation and write $V(v,z)u$ as $v^zu$. 
Let us see what several standard formulas look like in this new notation:
$$\eqalign{
\hbox{Old notation} &\qquad  \hbox{New notation}\cr
V(u,z)v &\qquad  u^zv\cr
V(V(a,x)b,y)c=V(a,x+y)V(b,y)c&\qquad (a^xb)^yc=a^{xy}(b^yc)\cr
V(a,x)V(b,y)c=V(b,y)V(a,x)c &\qquad  a^xb^yc=b^ya^xc\cr
V(a,x)b=e^{xL_{-1}}(V(b,-x)a)&\qquad  a^xb=(b^{x^{-1}}a)^x\cr
V(1,x)b=b &\qquad  1^xb=b\cr
}$$
The formulas on the right hand side are all easy to recognize: they 
are just standard formulas for a commutative ring acted on by a group $G$,
where $a,b,c$ are in the ring $V$, $x,y,z$ are elements of the group $G$,
and the action of $x\in G$ on $a\in V$ is denoted by $a^x$. This suggests
that we should try to set things up so that vertex algebras 
are exactly the commutative rings objects over some sort of 
mysterious group-like thing $G$. 

For simplicity we will work over a field of characteristic 0. This is
not an important assumption; it just saves us from some minor
technicalities about divided powers of derivations.

We will first look at the special case of vertex algebras such that
all the vertex operators $V(a,x)$ are holomorphic. We show that such
vertex algebras are the same as commutative algebras with a derivation
$D$. The correspondence is given as follows. First suppose that $V$ is
a commutative algebra with derivation $D$. We define the vertex
operator $V(a,x)$ by $V(a,x)b= \sum_{i\ge 0} (D^ia)bx^i/i!$.
Conversely if $V$ is a vertex algebra we define the product by
$ab=V(a,0)b$ and the derivation by $Da=$ coefficient of $x^1$ in
$V(a,x)b$.  (We cannot really check that this turns commutative
algebras into holomorphic vertex algebras and vice versa because we
have not yet said exactly what the axioms for a vertex algebra are.)

In the new notation for vertex algebras above we would put
$a^x=\sum_ix^iD^ia/i!$, $a^xb=\sum_ix^iD^iab/i!$. Here we think of $x$
as being an ``element'' of the one dimensional formal group $\hat
G_a$. This formal group has as its formal group ring $H$ the algebra
of polynomials $k[D]$ and its coordinate ring is the ring of formal
power series $k[[x]]$.  (In characteristic 0 it does not matter
whether we use Lie algebras or formal groups which are essentially
equivalent, but in other characteristics formal groups are better than
Lie algebras.) The (tensor) category of modules with a derivation is
the same as the category of modules over the formal group ring $H$, so
holomorphic vertex algebras are the same as the commutative ring
objects in this category.

What is the difference between a commutative algebra over $\hat G_a$
and a (non holomorphic) vertex algebra? The only difference is that
expressions like $a^xb^yc$ are no longer holomorphic in $x,y$ but can
have singularities; more precisely $a^xb^yc$ lies in
$V[[x,y]][x^{-1},y^{-1},(x-y)^{-1}]$.  In other words we can
provisionally define a vertex algebra to be a module $V$ such that we
are given functions $a^xb^yc^z\cdots$ for each $a,b,c,\ldots\in V$
which behave just like the corresponding functions for commutative
rings over $\hat G_a$, except that they are allowed to have certain
sorts of singularities. Notice that we can no longer reconstruct
a commutative ring structure on $V$ by defining $ab=a^xb$ at $x=1$,
because $a^xb$ may have a singularity at $x=1$. 

The definition above is too vague to be useful, so we try to make it
more precise. What we really want to do is to define some sort of
category, whose multilinear maps are somehow allowed to have the sort
of singularities above, and whose commutative ring objects are just
vertex algebras. We first ask in what sort of categories we can define
commutative ring objects. The obvious answer is tensor categories,
such as the category of modules over $\hat G_a$ (or over any
cocommutative Hopf algebra) but this turns out to be too
restrictive. (We will implicitly assume that all categories are
additive and have some sort of ``symmetric'' structure.)  A tensor
category requires that multilinear maps should be representable, but
this is sometimes not the case for the categories we are interested
in, and in any case this assumption is unnecessary. It is sufficient
to assume that for each collection of objects $A_1,\ldots, A_n,B$ of
the category we are given the space of multilinear maps
$Multi(A_1,\ldots, A_n,B)$, and that these satisfy a large number of
fairly obvious properties which I cannot be bothered to write down.

Unfortunately multilinear categories are not really the right objects
either.  The problem is the following. An expression like $a^xb^yc$
should live in a space like
$V[[x,y]][x^{-1},y^{-1},(x-y)^{-1}]$. However the expression
$a^x(b^yc)$ does not naturally live in this space, but in the
larger space $V[[y]][y^{-1}][[x]][x^{-1}]$, and $(a^{xy^{-1}}b)^yc$
lives in a different larger space. This makes it hard to compare these
expressions in a clean way. The easiest way to solve this problem is to
define it out of existence by using ``relaxed multilinear
categories''. The key idea is that instead of just once space of
multilinear maps $Multi(A_1,\ldots, A_n,B)$ we are given many
different spaces $Multi_p(A_1,\ldots, A_n,B)$ of multilinear maps,
parameterized by trees $p$ with a root (corresponding to $B$) and $n$
leaves, corresponding to $A_1,\ldots, A_n$. We should also have some
extra structures, consisting of maps between different spaces of
multilinear maps corresponding to collapsing maps between trees, and a
composition of multilinear maps taking multilinear maps of types
$p_1,\ldots, p_n, p$ to a multilinear map of type $p(p_1,\ldots,
p_n)$. (Here $p(p_1,\ldots,p_n)$ is 
 the tree obtained by attaching $p_1,\ldots, p_n$ to
the leaves of $p$.)  For details see my paper or Soibelman's notes, or
better still work them out for yourself.

The main point is that in a relaxed multilinear category it is still
possible to define commutative associative algebras. Joyal has pointed
out that the definition of associative algebras in a relaxed
multilinear category is strikingly similar to the definition of an
$A_\infty$ algebra; for example, the cells of the complexes used to
define $A_\infty$ algebras are parameterized by rooted trees with $n$
leaves, and the boundary maps correspond to the collapsing maps between trees. 

One way of constructing relaxed multilinear categories is as the
representations of ``vertex groups''. A vertex group can be thought of
informally as a group together with certain sorts of allowed
singularities of functions on the group. More precisely a vertex group
is given by a cocommutative Hopf algebra $H$, which we thing of as its
group ring, together with an algebra of ``singular functions'' $K$
over the ``coordinate ring'' $H^*$ of $H$.  The axioms for a vertex
group say that $K$ behaves as if it were the ring of meromorphic
functions over the ``group'' $G$; for example, the ring of meromorphic
functions is acted on by left and right translations, so $K$ should
have good left and right actions of $H$. A typical example of 
a vertex group is to take $H= k[D]$, $H^*=k[[x]]$ (so that $H$ is the formal group ring of $\hat G_a$), and to take $K$ to be the quotient field 
$k[[x]][x^{-1}]$ of $H^*$, which we can think of as the field of
rational functions on the formal group $\hat G_a$.

We can construct a relaxed multilinear category from a vertex group
roughly as follows. The underlying category is the same as that of the
Hopf algebra of the vertex group. However the spaces of multilinear
maps are different. Rather than define these in general, which is
a bit complicated, we will just look at one example. We take $G$ to
be the vertex group above (whose commutative rings are vertex
algebras), and take 3 $G$-modules $A$, $B$, and $C$. Then the space of
bilinear maps from $A,B$ to $C$ is defined to be the ordinary space of
bilinear maps from $A\times B$ to $C[[x,y]][(x-y)^{-1}]$ which are
invariant under an action of $G^3$. (The easiest way to work out what
the action of $G^3$ should be is to see what it has to be for the
invariant bilinear maps from $A\times B$ to $C[[x,y]]$ 
taking $a\times b$ to $\sum_{i,j}f(D^ia,D^jb)x^iy^j/i!j!$ to be the same
as invariant maps $f$ from $A\times B$ to $C$.)

We summarize what we have done so far:
\item 1 We have introduced ``vertex groups'''.
\item 2 The modules over a vertex group form a ``relaxed multilinear
category''.
\item 3 The commutative ring objects over the simplest nontrivial 
vertex group are exactly vertex algebras.

Now that we have set up this machinery, it is easy to find higher
dimensional analogues of vertex algebras: all we have to do is look at
commutative algebras over higher dimensional vertex groups $G$; we
will call these vertex $G$ algebras.  As an example we will construct
vertex algebras related to free quantum field theories in higher
dimensions.  (Can one construct vertex algebras corresponding to
nontrivial quantum field theories in higher dimensions? At the moment
this is just a daydream, as it is too vague to be called a
conjecture.)

We first need to construct a suitable vertex group $G$. We take its
underlying Hopf algebra $H$ to be the polynomial algebra
$\R[D_0,\ldots, D_n]$ where $D_i=\partial /\partial x_i$, which we
think of as the universal enveloping algebra of the Lie algebra of
translations of spacetime (with $x_0=t$). The dual $H^*$ is then
$\R[[x_0,\ldots, x_n]]$, which we think of as the algebra of functions
on spacetime. We define $K$ to be $H^*[(x_0^2-x_1^2-\cdots
x_n^2)^{-1}]$, which we think of as the algebra of functions on
spacetime which are allowed to have singularities (poles) on the light
cone.

Now we define the vertex $G$ algebra $V$. The underlying space of $V$
is the universal commutative $H$-algebra generated by an element
$\phi$, so $V=\R[\phi, D_0\phi,D_1\phi,\ldots, D_0^2\phi,\ldots]$ is a
ring of polynomials in an infinite number of variables. We think of
$V$ as the ring of classical fields generated by $\phi$, and it is a
(holomorphic) vertex $G$ algebra as it is a commutative ring acted on
by $H$. We will turn it into a a nontrivial vertex $G$ algebra by
``deforming'' this trivial vertex $G$ algebra structure. (In general,
for vertex $G$-algebras, quantization means deforming the structure on
some commutative ring to turn it into a vertex $G$ algebra.) 

To do this we recall the following method of constructing commutative
rings: if $V$ is a space acted on by commuting operators $v_n$, and if
$V$ is generated by an element $1\in V$ by the action of these
operators, then $V$ has a unique commutative ring structure such that
1 is the identity and the actions of all the operators are given by
multiplication by elements of $V$.  (Proof: easy exercise.) A similar
theorem holds for vertex algebras (as was proved by Frenkel, Kac,
Rado, and Wang). We will make $V$ into a vertex $G$ algebra by finding
a vertex operator $\phi(x)=\phi(x_0\ldots, x_n)$ acting on $V$ such
that $\phi(x)\phi(y)=\phi(y)\phi(x)$ and applying the construction
above.

To construct $\phi(x)$, we first put
$$\phi^+(x)=\sum_iD^i\phi x^i/i!.$$
The vertex algebra structure on $V$ defined by this vertex operator is just
the commutative ring structure on $V$, so we need to deform $\phi^+$. 
We define $\phi^-(x)$ to be the unique $G$-invariant derivation 
from $V$ to $V[x][(x_0^2-x_1^2\cdots)^{-1}]$ taking $\phi$ to some
even function $\Delta(x)$ (called the propagator). This is uniquely defined 
by the universal property of $V$. Finally we put
$$\phi(x)=\phi^+(x)+\phi^-(x).$$  It is easy to check that $\phi(x)$
and $\phi(y)$ commute, as $\phi^+(x)$ and $\phi^+(y)$ commute,
$\phi^-(x)$ and $\phi^-(y)$ commute, and
$[\phi^-(x),\phi^+(y)]=-[\phi^+(x),\phi^-(y)]=
\Delta(x-y)$. Therefore we can make $V$ into a
commutative vertex $G$-algebra.

Notice that in quantum field theory, $\phi(x)$ means the value of some
operator valued distribution $\phi$ at some point $x$ of a manifold.
On the other hand, for vertex $G$ algebras, $\phi(x)$ should be
thought of as the action of a ``group'' element $x$ on an element
$\phi$ of a vertex $G$ algebra. Several other concepts in quantum
field theory can also be translated into vertex algebra theory; for
example, the correlation functions are
$Tr(\phi(x)\phi(y)\phi(z)\cdots)$, where $Tr$ is some $G$ invariant
linear function on $V$. Some concepts are not so easy to extend; for
example, vertex $G$ algebras do not seem to be able to cope with
arbitrary curved spacetimes other than Lie groups, and it can be
difficult to reconstruct a Hilbert space (this includes as a very
special case the problem of deciding which representations of the
Virasoro algebra are unitary).

Finally we will briefly describe a special case of a theorem about vertex
$G$-algebras generalizing the usual identity for vertex algebras.
This theorem says (roughly) that under fairly general conditions, the
integral of a vertex operator over an $n$-cycle is a vertex
differential operator of order $n$. We will illustrate this for
ordinary vertex algebras, when it says that the integral of $a(x)$
around the origin is a vertex differential operator of order 1.
This means 
$$\int a(x)dxb(y)-b(y)\int a(x)dx$$ 
has order 0 and is therefore of the form $t(y)$ for some $t$.
Applying both sides to 1 and taking $y=0$ shows that $t=\int a(x)dxb$. 
Therefore
$$\int a(x)dxb(y)c-b(y)\int a(x)dxc= (\int a(x)dxb)(y) c.$$ 
Finally integrating both sides around $y=0$ shows that
$$\int a(x)dx\int b(y)dyc-\int b(y)dy\int a(x)dxc= \int (\int
a(x)dxb)(y)dy c.$$ which is a special case of the vertex algebra
identity. (The other cases can be deduced in a similar way.)

\bye